\documentclass[12pt]{amsart}
\usepackage{amscd,amsthm,amssymb,amsfonts,amsmath,epsfig,epic,eclbip}
\usepackage{bbm}
\usepackage{url}
\usepackage{color}
\usepackage[arrow,matrix,tips,2cell]{xy}

\theoremstyle{plain}
\newtheorem{thm}{Theorem}[section]
\newtheorem{lemma}[thm]{Lemma}
\newtheorem{prop}[thm]{Proposition}

\newtheorem{conj}[thm]{Conjecture}
\theoremstyle{definition}
\newtheorem*{defn}{Definition}

\theoremstyle{remark}
\newtheorem*{remark}{Remark}
\newtheorem*{ack}{Acknowledgments}
\newtheorem*{notations}{Notations}

\newcommand{\Z}{\mathbb Z}    
\newcommand{\R}{\mathbb R}    
\newcommand{\C}{\mathbb C}    
\newcommand{\A}{\mathbb A}    
\newcommand{\Af}{\mathcal Af\negmedspace f}
\newcommand{\SL}{\operatorname{SL}}
\newcommand{\GL}{\operatorname{GL}}
\newcommand{\Aff}{\operatorname{Af{}f}}
\newcommand{\Arg}{\operatorname{Arg}}

\newcommand{\Vol}{\operatorname{Vol}}
\newcommand{\Hess}{\operatorname{Hess}}

\newcommand{\<}{\langle}   
\renewcommand{\>}{\rangle} 

\newcommand{\D}{\Delta}
\newcommand{\DSv}{\Delta^\vee_\lambda}
\newcommand{\Dv}{\Delta^\vee}
\newcommand{\DT}{{\Delta_\nu}}
\newcommand{\smooth}{{\Sigma\backslash D}}
\newcommand{\carrier}{\operatorname{carrier}}

\newcommand{\suchthat}{\ : \ }

\newcommand{\dD}{\partial\Delta}
\newcommand{\dDv}{\partial\Delta^\vee}
\newcommand{\U}{\mathcal U}
\newcommand{\V}{\mathcal V}
\newcommand{\T}{\mathbb T}

\newcommand{\am}{\mathcal A^\lambda}

\newcommand{\aff}{{\operatorname{af{}f}}}
\newcommand{\join}{\ast}
\renewcommand{\setminus}{\backslash}

\newcommand{\vertS}{\operatorname{vert}(S)}
\newcommand{\vertT}{\operatorname{vert}(T)}

\newcommand{\dDT}{\partial \DT}
\renewcommand{\L}{\mathcal L}

\newcommand{\Hsm}{Z_a^{\mathrm{sm}}}

\newcommand{\dDS}{\partial \DSv} 
\newcommand{\dU}{\partial\U}
\newcommand{\dV}{\partial\V}

\newcommand{\Cstar}{\C\setminus\{0\}}
\newcommand{\cone}{\operatorname{cone}}

\newcommand{\F}{\mathcal F}
\newcommand{\X}{\mathfrak X}

\newcommand{\SC}{\operatorname{SC}}
\newcommand{\sm}{\mathrm{sm}}

\renewcommand{\caption}[1]{%
  \refstepcounter{figure}
  \begin{minipage}[t]{.8\textwidth}
    \vspace{.3\baselineskip}
    \begin{center}
      \scriptsize {\sc Figure \thefigure:}
      #1
    \end{center}
  \end{minipage}
}

\begin{document}

\begin{abstract}
  This paper is a continuation of our paper \cite{HZh} where we have
  built a combinatorial model for the torus fibrations of Calabi-Yau
  toric hypersurfaces. This part addresses the connection between the
  model torus fibration and the complex and K\"ahler geometry of the
  hypersurfaces.
\end{abstract}

\title[Integral affine structures on spheres II]{Integral affine
  structures on spheres 
\makebox[0pt]{\raisebox{1in}{\normalfont DUKE-CGTP-03-01}} \\
and torus fibrations of Calabi-Yau toric
  hypersurfaces II} 
\author{Christian Haase} 
\author{Ilia Zharkov}
\address{Mathematics Department \\ Duke University \\ Durham, NC 27708
  \\ USA}
\email{[haase,zharkov]@math.duke.edu}
\thanks{The first author gratefully acknoledges support by NSF-grant
  DMS-0200740.}
\maketitle

\section{Introduction}
In this paper we endow the topological model $\Sigma$ constructed in
\cite{HZh} with the structure of a metric space which is a K\"ahler
affine manifold away from the (codimension 2) discriminant locus and
relate it to the {\it geometry} of the Calabi-Yau family of toric
hypersurfaces $Z_s$ near large complex structure limit point.

The main result is the following. Given an integral K\"ahler affine
structure on $\smooth$ which is in the right class and satisfies
certain (bi-polyhedral) compatibility conditions, we construct a
family of K\"ahler metrics on $Z_s$ such that as $s\to\infty$ (the
large complex structure limit):
\begin{itemize}
\item The embedding $\psi: Z_s^\sm\hookrightarrow W_s$ of ``smooth''
  portions of the hypersurfaces into the model torus bundle (with the
  right twist) identifies the scalar products and the complex
  structures on the tangent spaces $T_x Z_s$ and $T_{\phi(x)} W_s$ up
  to terms of order $o(1)$ uniformly in $x\in Z_s^\sm$.
\item the Gromov-Hausdorff distance between the pairs
  $(Z_s,Z_s\setminus Z_s^\sm)$ and $(\Sigma,D)$ is of order $o(1)$.
\end{itemize}

Interchanging the input data $(\D,S,\lambda)$ with $(\Dv,T,\nu)$ and
repeating the construction of the K\"ahler metrics for the {\em dual}
Calabi-Yau family gives rise to the same limiting metric space
$\Sigma$ and the {\em dual} K\"ahler affine structure.
This result constitutes a significant part of the the limiting mirror
symmetry conjecture (cf. \cite{KS}).

The major missing part toward proving the metric collapse is that our
metrics are not Ricci-flat. To establish the Ricci-flatness away from
the discriminant one would need the affine Calabi conjecture (see
Section~\ref{sec:calabi}). But even assuming a Monge-Amp\`ere solution
on $\smooth$ the behavior of true Calabi-Yau metrics is not expected
to be approximated by a semi-flat construction near the discriminant
locus. A further development in this direction requires some local
estimates on CY metrics near singularities. We hope to address these
issues by using generalized Gibbons-Hawking ansatz elsewhere
(cf. \cite{Zh2}).

\begin{notations}
We continue to use notations from \cite{HZh}.
\begin{itemize}
  \item $(\D,S),(\Dv,T)$ is a dual pair of $d$-dimensional reflexive
    polytopes with coherent triangulations of their boundaries. 
  \item $\D_\Z,\Dv_\Z$ are the sets of integral points of $\D,\Dv$.
  \item $\SC(S),\SC(T)$ are the secondary cones in $\R^{\D_\Z},\R^{\Dv_\Z}$    
    corresponding to the triangulations $S\join\{0\},T\join\{0\}$ of $\D,\Dv$. 
  \item $\lambda$ and $\nu$ are integral vectors in the interiors of
    the respective secondary cones.
  \item $X_T$ is the toric variety associated to the simplicial fan in
    $(\R^d)^*$ given by the triangulation $T$.
\end{itemize}
\end{notations}

\begin{ack}
  We would like to thank Lev Borisov, Robert Bryant, Mark Gross, Maxim
  Kontsevich, Grisha Mikhalkin, Dave Morrison, Svetlana Roudenko, Mark
  Stern and Stephanos Venakides for valuable conversations, and the
  entire Duke CGTP group for a stimulating environment.
\end{ack}

\section{The model: K\"ahler affine structures and torus bundles}
\subsection{Integral K\"ahler affine structures}
Let $\A^n$ be $n$-dimensional affine space.
An {\it integral affine structure} on an $n$-dimensional manifold $Y$
is given by an open covering $\{U_\alpha\}$ of $Y$ together with
coordinates $\phi_\alpha: U_\alpha\to\A^n$ such that the transition
maps $\phi_\alpha\circ \phi_\beta^{-1}$ are in $\SL(n,\Z) \ltimes
\R^n$ on the non-empty overlaps $U_\alpha \cap U_\beta$.
An integral {\it K\"ahler} affine structure on $Y$ is a Riemannian
metric $g$ which is potential in local affine coordinates, i.e.
$g_{ij}=\frac{\partial^2K_\alpha}{\partial y_i \partial y_j}$ for some
local potentials $K_\alpha$.

The dual K\"ahler affine structure on the same Riemannian manifold
$(Y,g)$ is defined as follows (cf \cite{KS}). We use the same covering
$\{U_\alpha\}$. The new affine coordinates are $\hat y_i=\frac{\partial
  K_\alpha}{\partial y_i}$ which take values in the dual affine space
$({\A^n})^*$ (the underlying vector spaces for ${\A^n}$ and $({\A^n})^*$
are naturally dual). The new local potentials $\hat{K}_\alpha$ are
defined by the Legendre transforms of the old ones: 
\begin{equation*}
 \hat{K}_\alpha(\hat y)=\max_{y\in U_\alpha}\{\<\hat y,y\>-K_\alpha(y)\}. 
\end{equation*} 
Here one needs to choose origins in ${\A^n}$ and $(\A^n)^*$ to define
the pairing. Different choices give rise to equivalent K\"ahler affine
structures. The dual affine structure is integral iff the original one is.
 
Given an affine structure on $Y$ one can consider its monodromy
representation $\pi_1(Y)\to SL(n,\Z) \ltimes \R^n$.  Two equivalent
affine structures have conjugate monodromies. 

For a K\"ahler affine manifold $(Y,g)$ one can define (cf. \cite{KS})
a characteristic class $[g]$ of the metric, which is an analog of the
K\"ahler class in complex geometry. Let $\Af_Y$ be the sheaf of
locally affine functions. The metric is given by local potentials in
affine coordinates: $g_{ij}=\frac{\partial^2K}{\partial y_i \partial
  y_j}$. Then the differences of the potentials on the overlaps will
define a \v Cech cohomology class $[g]\in H^1(Y, \Af_Y)$.

It is more natural to combine the monodromy representation and the
metric class into one class, which we will call the class of {\it
  affine polarization}. It can be represented by a \v{C}ech 1-cocycle
with values in the semi-direct product
$(\SL(n,\Z)\ltimes\R^n)\ltimes\Aff_n$, where the affine
transformations $\SL(n,\Z) \ltimes \R^n$ act on the affine functions
$\Aff_n$ from the right.

The natural projection onto the normal component
$\SL(n,\Z)\ltimes\R^n$ in the above semi-direct product gives the
monodromy representation. To recover the metric class, however, one
needs to fix a splitting of the natural map
$(\SL(n,\Z)\ltimes\R^n)\ltimes\Aff_n\to\Aff_n$. Different splittings
will give conjugate metric classes.

For the purposes of this paper we consider a convenient
$(n+2)$-dimensional faithful representation of the group
$(\SL(n,\Z)\ltimes\R^n)\ltimes\Aff_n$. Let us choose $q$ -- an
integral vector in $\R^{n+2}$, and $p$ -- an integral vector in the
dual space $(\R^{n+2})^*$. Then this representation provides an
isomorphism of $(\SL(n,\Z)\ltimes\R^n)\ltimes\Aff_n$ with the
following subgroup of $\GL_{n+2}(\R)$:
\begin{equation*} 
  G_n(p,q):=\{g\in\GL_{n+2}(\R)\suchthat g(q)=q,\ g^*(p)=p,\
  g|_{\{\<p,x\>=0\}/q} \text{ is integral}\},
\end{equation*}
where $g^*:(\R^{n+2})^*\to(\R^{n+2})^*$ is the adjoint linear
transformation. The affine space $\A^n$ can be identified with
$\{\<p,x\>=1\}/q$, and the $G_n(p,q)$ action on it gives the corresponding
affine transformation of $\A^n$. To recover the affine function
$f:\A^n\to\R$ one needs to fix an integral linear functional
$l\in(\Z^{n+2})^*$, such that $l(q)=1$. Then $f(x)=l(g(x))-l(x)$ is a
function, well defined on the quotient $\{\<v,x\>=1\}/w$. Only the
representing \v{C}ech cocycle depends on the choice of $l$, not the
metric class itself.

For the (mirror) symmetry sake we also choose an integral element
$k\in \R^{n+2}$ with $p(k)=1$. The vector $k$ defines an origin in
$\A^n$, hence it allows to recover the translational part of the
affine transformation. Again, the class of this translational part,
called the {\it radiance obstruction} (cf. \cite{GoldHir}), is
independent of the choice of $k$ (different $k$'s give rise to
conjugate monodromies).
As was noted in \cite{GS} the radiance obstruction class is dual to
the linear part of the metric class under the duality between the
K\"ahler affine structures.

More generally, it is also clear that the full polarization class for
the dual K\"ahler affine structure can be represented by the adjoint
inverse transformations for each $U_\alpha \cap U_\beta$, with the
r\^oles of $q,k$ and $p,l$ exchanged.
Given a basis $\{e_i\}$ of $\R^{n+2}$ such that $\<p,e_i\>=0,\
i=1,\dots,n+1$, $e_{n+1}=q$ and $e_{n+2}=k$, the group $G_n(p,q)$ can
be represented by non-degenerate matrices in the form
\begin{equation*}
\left(
\begin{array}{ccc}
A & 0 & b \\ 
a & 1 & c \\
0 & 0 & 1
\end{array}
\right),
\end{equation*}
where $A$ and $b$ represent the linear and translational parts of the
affine transformations, and $a,c$ are the linear and constant parts of
the affine function, respectively.

All of the above (including the affine structure itself) can be
defined even if we do not require the affine charts to be maps into
the {\it same} affine space. We won't have groups anymore, but in all
cocycle conditions the compositions still make sense. In particular,
to specify an integral affine structure we would need continuous maps
$\phi_\alpha: U_\alpha\to\A^n_\alpha\cong
\{\<p_\alpha,x_\alpha\>=1\}/q_\alpha$ together with integral elements
$q_\alpha,k_\alpha\in \R^{n+2}_\alpha$ and $p_\alpha,l_\alpha\in
(\R^{n+2})^*$ with $\ p_\alpha(k_\alpha)=1,\ l_\alpha(q_\alpha)=1$,
such that the transition maps
$\phi_{\alpha\beta}:{\R^{n+2}}\to{\R^{n+2}}$ satisfy the corresponding
invariance, coinvariance and integrality conditions.

The cohomological information, such as monodromy, radiance obstruction
and the metric class, is encoded in the transition maps
$\phi_{\alpha\beta}$.

\subsection{Bi-polyhedral K\"ahler affine structures}\label{sec:bipikas}
We will be interested in a very special types of integral K\"ahler
affine structures. These structures arise in the metric limits of
Calabi-Yau hypersurfaces and complete intersections in toric
varieties.
\begin{defn}
  An integral affine structure on $Y$ is {\it polyhedral} if there is an
  $n$-di\-men\-sio\-nal polyhedral complex $P$, a collection of disjoint
  open sets $\{U_\alpha\}$, whose closures cover $Y$, i.e. $Y=\bigcup
  \overline{U}_\alpha$, and a continuous map $\phi:Y\to P$, which
  provides an affine homeomorphism of each $U_\alpha$ with the
  interior of some $n$-dimensional face of $P$. We say that the pair
  $(\{U_\alpha\},P)$ realizes the polyhedral affine structure if $P$
  is minimal, which, in particular, means that there is a bijection
  between open sets $\{U_\alpha\}$ and $n$-dimensional cells of $P$.
\end{defn}

\begin{defn}
  An integral K\"ahler affine structure on $Y$ is {\it bi-polyhedral}
  (bi-PIKAS for short)
  if there is a bipartite covering $\{U_\alpha,V_\beta\}$ of $Y$ and
  two polyhedral complexes $P,\hat{P}$ such that $(\{U_\alpha\},P)$
  and $(\{V_\beta\},\hat P)$ provide polyhedral realizations of the
  underlying affine structure and its dual, respectively. We say that
  the bi-polyhedral K\"ahler affine structure is of type
  $(\{U_\alpha,V_\beta\},P,\hat P)$.
\end{defn}

The bi-polyhedral property imposes very severe restrictions on the
compatibility between Riemannian metric and affine structure. In
particular, $\Hess K_\alpha\in\L^1_{loc}(\R^n)$  and the
Cauchy-Schwartz inequality implies that the
metric completion of $Y$ can be identified with $P$ or $\hat P$. This
endows both polyhedral complexes with (isomorphic) structures of complete
metric spaces.

Next we want to show the existence of bi-PIKAS on $\smooth$.
Recall from \cite{HZh} that $\smooth$ has a bipartite covering by open
sets $U_v$ and $V_w$. Also, given vectors $\lambda,\nu$ in the
interiors of the respective secondary cones $\SC(S),\SC(T)$ with
$\lambda(0)=\nu(0)=0$ we can define the polytopes
\begin{gather*}
  \DSv = \{ n \in \R^d \suchthat \<m,n\> + \lambda(m) \le 0
  \text{ for all } m \in \D_\Z\}, \\
  \DT = \{ m \in (\R^d)^* \suchthat \<m,n\> + \nu(n) \le 0
  \text{ for all } n \in \Dv_\Z\}.
\end{gather*}
In the future we will abbreviate the type of a bi-polyhedral integral
K\"ahler affine structure on $\smooth$ by simply $(\lambda,\nu)$
having fixed the covering $(\{U_v,V_w\})$.

In order to specify a bi-PIKAS of type $(\lambda,\nu)$ on
$\smooth$ we will provide
the following data. A Legendre dual pair of convex functions
$\Phi,\hat\Phi$ on $\DSv,\DT$, respectively, smooth on each strata of
the respective polytope. (This implies that the Hessians of both
$\Phi,\hat\Phi$ are positive along the strata.) Then the restrictions
of $\Phi$ to the facets of $\DSv$ serve as potentials for the metric
along $U_v$. 
For future use we will prove
that it is possible to choose these functions consistently in
$\lambda$ and $\nu$.
\begin{defn}
  Suppose, for each pair $(\lambda,\nu)\in\SC(S)\times\SC(T)$ we have
  a bi-polyhedral K\"ahler affine structure on $\smooth$ of type
  $(\lambda,\nu)$, which varies continuously with $(\lambda,\nu)$ in
  the Hausdorff topology of metric structures on $\Sigma$. We call
  such a family {\it projective} if:
  \begin{itemize} 
  \item For any linear functions $l\in\R^d,\ l^\vee\in(\R^d)^*$,
    the bi-PIKAS for
    $(\lambda+l,\nu+l^\vee)\in\SC(S)\times\SC(T)$ have the same
    underlying K\"ahler affine structure.
  \item The bi-PIKAS for $(\epsilon^{-1}\lambda,\epsilon\nu)$ differs
    from the bi-PIKAS for $(\lambda,\nu)$ by the $\epsilon$-rescaling
    \begin{equation*}
(y'_\alpha,\hat y'_\alpha, K'_\alpha(y'_\alpha), \hat K'_\alpha(\hat y'_\alpha),g_{ij})= (\epsilon^{-1}y_\alpha,\epsilon \hat y_\alpha, K_\alpha(\epsilon y'_\alpha), \hat K_\alpha(\epsilon^{-1}\hat y'_\alpha),\epsilon ^2 g_{ij})
\end{equation*}
  \end{itemize}
\end{defn}
Note here that adding global linear functions to the potentials
$\Phi,\hat\Phi$ will induce translations of the polytopes
$\DT,\DSv$. Though giving different bi-PIKAS (as we defined them) this
will have no effect on the underlying K\"ahler affine structures (the
latter will be canonically equivalent).

Another important observation is that rescaling the data for bi-PIKAS will provide the same metric on
$\smooth$, though different affine structures.
\begin{prop} \label{prop:existence}
  There are Legendre dual functions $\Phi \colon \DSv \rightarrow \R$
  and $\hat\Phi \colon \DT \rightarrow \R$ that define a projective
  family of bi-PIKAS on $\smooth$.
\end{prop}
\enlargethispage{10mm}
\begin{proof}
  First, we choose a smooth function with positive Hessian on
  $\Delta_{\lambda-\beta}^\vee$ whose gradients stay in
  $\Delta_{\nu-\beta}$, and cover $\Delta_{\nu-2\beta}$.
  \begin{center}
    \includegraphics{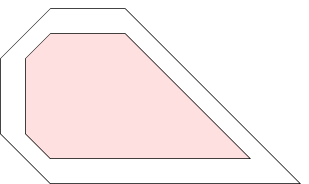} \qquad \qquad
    \includegraphics{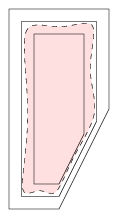} \\
    \caption{\quad The domain for the first step. \qquad The set of
      gradients.}
  \end{center}
  \pagebreak
  As a second step, we need a continuous strictly convex function
  $\tilde{\nu}$ on $\beta \Dv$ that is an approximation of the
  function $\nu$ with the following properties.
  \begin{itemize}
    \item $\tilde{\nu}$ is piecewise smooth.
    \item $\operatorname{Hess} \tilde\nu > 0$ on the smooth pieces.
    \item The gradients along $\cone \tau$ belong to a neighborhood of
    the corresponding vertex in $\DT$ that are pairwise disjoint, and
    do not meet the $\beta$ neighborhood of the barycenter of $\DT$.
    \item For a vertex $w \in T$, the $w$-directional derivatives
    equal $-\nu(w)$ in the star neighborhood of $\beta w$ in the
    barycentric subdivision of $\beta T$.
  \end{itemize}
  \begin{center}
    \includegraphics{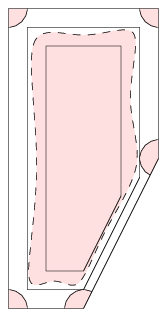} \\
    \caption{The set of gradients of $\tilde\nu$.}
  \end{center}
  We obtain a convex function on $(\R^d)^*$ if we consider the
  lower hull of the ($d+1$)-dimensional Minkowski sum of graphs of the
  two functions.

  Finally, we want to obtain a function that is smooth along the
  strata. As explained in \S~\ref{sec:regularization}, we convolute
  with a kernel that depends on the position $x \in \DSv$ as follows.
  Consider the quadratic form
  \begin{equation*}
    Q(x) = \sum_{v \in S} Q_v(x) \ \text{ where } \ Q_v(x) =
    \frac{1}{\rho_v(x)} v^2
  \end{equation*}
  This quadratic form is non-degenerate because the $v$'s span
  $(\R^d)^*$. It has a dominant summand if $x$ is close to a facet.
  Now our kernel will be a normalized $e^{1-Q(x)}$. Its support -- the
  ellipsoid given by $Q(x) \le 1$ -- depends on the position as
  sketched in the figure.
  \begin{center}
    \includegraphics{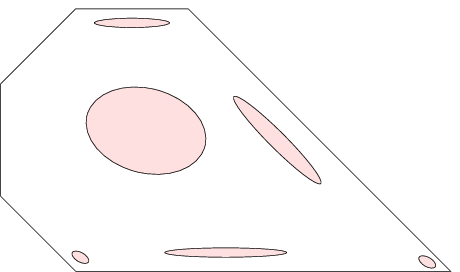} \\
    \caption{The support of the mollifier.}
  \end{center}
  The obtained function will have a positive Hessian along the
  strata, so that the Legendre dual function will be smooth along its
  corresponding strata.
\end{proof}

The metric completion of $\smooth$ can be identified with $\Sigma$ and
endowed with the structure of a compact metric space via the
bi-polyhedral homeomorphisms $\phi,\hat\phi$:
\begin{equation*} 
  \xymatrix{
    & \Sigma \ar[dl]_{\phi} \ar[dr]^{\hat\phi}&\\  
    \dDS & & \dDT
  }
\end{equation*}
Throughout the paper we will often identify points in $\Sigma$, $\dDS$
and $\dDT$ by means of these homeomorphisms when there is no
confusion.

We can realize the affine coordinates $y$ on $U_v$ and $V_w$
explicitly as taking values in the following (affine) subspaces and
quotients of $\R^d$:
\begin{gather*}
  y(q)\in \R^d_v(\lambda(v)):=\{ n\in\R^d \suchthat
  \<v,n\>+\lambda(v)=0\},\quad q\in U_v,\\
  y(q)\in \R^d/w, \quad q\in V_w,
\end{gather*}
and the transition maps are given by the obvious projections
$\R^d_v(\lambda(v))\to\R^d/w$. Then the affine
monodromy along a primary loop $(v_0w_0v_1w_1)$ is given by
(cf. \cite[Lemma 2.4]{HZh}):
\begin{equation}\label{eq:monodromy}
  n\mapsto n+[\<v_1,n\>+\lambda(v_1)](w_1-w_0), \quad n\in
  \R^d_v(\lambda(v)).
\end{equation}

To describe the full polarization class of a
bi-PIKAS of type
$(\lambda,\nu)$ on $\smooth$ we consider the representation of
$(\SL(d-1,\Z)\ltimes\R^{d-1})\ltimes\Aff_{d-1}$ in $\R^d\oplus\R$. The
dual space is identified with $(\R^d)^*\oplus\R$. For the charts $U_v$
and $V_w$ we set
\begin{gather*}
  U_v:\ p=\left(\begin{array}{c} 0\\1 \end{array} \right),\ 
  k=\left(\begin{array}{c} k_v\\0 \end{array} \right),\ 
  q=(v,0),\ l=(0,1),\\
  V_w:\ p=\left(\begin{array}{c} w\\0 \end{array} \right),\ 
  k=\left(\begin{array}{c} 0\\1 \end{array} \right),\ 
  q=(0,1),\ l=(l_w,0),
\end{gather*}
where we have chosen integral elements $l_w\in (\Z^d)^*$ and
$k_v\in\Z^d$ such that $\<l_w,w\>=1$ and $\<v,k_v\>=1$.
Then for $U_v\cap V_w$ the cocycle transformation $g_{vw}:
\R^d\oplus\R\to \R^d\oplus\R$ is given by:
\begin{equation*}
  \left(
    \begin{array}{c}
      n \\ s
    \end{array}
  \right)\mapsto \left(
    \begin{array}{c}
      n' - \left[1+\nu(w)\right] \<l_w,n'\>w+ sw\\ \<v,n\>
    \end{array}
  \right),
\end{equation*} 
where $n'=n-\<v,n\>\left[1+\lambda(v)\right]k_v$. The cocycle
$\{g_{vw}\}$ represents the polarization class which we will denote by
$[\lambda,\nu]$.

Then the monodromy representation $\pi_1(\smooth)\to
\GL(\R^d\oplus\R)$ along a primary loop $(v_0w_0v_1w_1)$ is given by
\begin{equation*}
  \left(
    \begin{array}{c}
      n \\ s
    \end{array}
  \right)\mapsto \left(
    \begin{array}{c}
      n \\ s
    \end{array}
  \right)
  +\alpha(n)\left(
    \begin{array}{c}
      w_1-w_0\\ \nu(w_1)-\nu(w_0)
    \end{array}
  \right),
\end{equation*} 
where
$\alpha(n)=v_1(n)+\lambda(v_1)v_0(n)-[1+\lambda(v)]v_1(k_{v_0})v_0(n)$.
Considering this transformation on the quotient by the last coordinate and using $k_{v_0}$
to identify $\A^{d-1}=\{\<v,n\>=1\}$ with $\R^d_v(\lambda(v))=\{\<v,n\>+\lambda(v)=0\}$ via
\begin{gather*}
n\mapsto n-\<v,n\>\left[1+\lambda(v)\right]k_v,
\end{gather*}
we recover the above affine monodromy on
$\R^d_v(\lambda(v))$.

Following through the above calculation shows that the converse is
also true:
any bi-polyhedral K\"ahler affine structure on $\smooth$ in the class
$[\lambda,\nu]$ is, in fact, of type $(\lambda,\nu)$.

\subsection{The Calabi conjecture}\label{sec:calabi}
Among all bi-polyhedral K\"ahler affine structures of type
$(\lambda,\nu)$ we expect to find a unique distinguished
representative -- the Monge-Amp\`ere structure: in affine coordinates
the metric satisfies $\det g_{ij}=c$. Its metric completion to
$\Sigma$ is supposed to be the limit of the Ricci-flat metrics on the
families of Calabi-Yau hypersurfaces.

\begin{conj}
[cf. also \cite{KTch}]
There is a unique bi-polyhedral K\"ahler affine structure on $\smooth$
of type  $(\lambda,\nu)$ such that the metric is Monge-Amp\`ere: $\det
g_{ij}= (\Vol \dDS)^{-1}\cdot\Vol\dDT$.
\end{conj}

Note that the Monge-Amp\`ere constant $c$ is determined from
calculating the metric volume of $\Sigma$ as $\sqrt{c}\cdot\Vol
\dDS=\sqrt{c^{-1}}\cdot\Vol\dDT$, where $\Vol$ means the affine volume
of the corresponding polytopal complex. Also note that the rescaled
data $(\epsilon^{-1}\lambda,\epsilon\nu,K_\alpha(\epsilon
y_\alpha),\hat K_\alpha(\epsilon^{-1}\hat y_\alpha), \epsilon^2 g_{ij})$
provides the Monge-Amp\`ere structure in the class
$[\epsilon^{-1}\lambda,\epsilon\nu]$ with the same metric on
$\Sigma$. Thus the Monge-Amp\`ere bi-PIKAS fit together into a
projective family.

We will not address this conjecture any further here, rather we will
be happy to start with any bi-polyhedral K\"ahler affine structure on
$\smooth$ given by a pair $(\Phi,\hat\Phi)$.

\subsection{The model torus fibrations as K\"ahler manifolds} 

The $(d-1)$-torus fibration $W(\lambda,\theta)$ (more naturally, a
$\T^{d-1}$-torsor) will depend on additional phase multi-parameter
$\theta:=\{\theta_v\}, v\in \vertS$, where all $\theta_v$ have values
in $\R/\Z$. First, we form (trivial) affine torus bundles over the
affine open sets by identifying the fibers with the affine tori
\begin{gather*}
  \T_v(\theta_v):=\{ n\in\T \suchthat \<v,n\>+\theta_v\equiv 0 \mod \Z
  \} \quad\text { over } U_v,\\
  \T/w:=(\R^d/w)/(\Z^d/w) \quad\text{ over } V_w.
\end{gather*}
The gluing maps are independent of a base point in an overlap $U_v\cap
V_w$ and given there by the natural projection $\T_v(\theta_v)\to
\T/w$. This defines the torsor $\pi: W(\lambda,\theta)\to\smooth$.

The topology of the total space of the torsor is determined by the
combinatorics of $\Sigma$, i.e., independent of $\lambda$ and $\theta$
as long as $\lambda$ is in the right secondary cone $\SC(S)$. In
particular, all $W(\lambda,\theta)$ are diffeomorphic to each other,
though not canonically.

Since the linear parts of the transition maps are the same for the
base and for the fibers, the tangent space at any point in
$W(\lambda,\theta)$  splits canonically as
$$T_{W(\lambda,\theta)}\cong T_\smooth\oplus T_\smooth.
$$  
This allows to define a canonical (integrable) almost complex
structure on $W(\lambda,\theta)$ as
\begin{equation*}
  J_w=\left( \begin{array}{cc}
      0 & \mathbbm{1} \\
      -\mathbbm{1} &0 \\ 
    \end{array} \right).
\end{equation*}

Given a Riemannian metric $g_{ij}$ on $Y$, one can define the pullback
metric $\pi^*(g_{ij})$ on $W(\lambda,\theta)$ which is, in fact,
K\"ahler. If, in addition, $g_{ij}$ satisfy the real Monge-Amp\`ere
equation, the induced metric on $W(\lambda,\theta)$ is Ricci-flat.

There is another slightly different description of the torus bundle
over a K\"ahler affine manifold $Y$ (cf. \cite{KS}), which is useful
when considering limiting behavior of Calabi-Yau degenerations.
We define a torus fibration $W_\epsilon(\lambda)$ as the quotient of
the total space of the tangent bundle $T({\smooth})$ by the integral
lattice spanned by $\left\{\epsilon\frac{\partial}{\partial
    y_i}\right\}$, where $y_i$ are the affine coordinates. This torus
bundle carries canonical complex structure and the pullback metric
which comes from the splitting $ T_{W_\epsilon(\lambda)}\cong
T_\smooth\oplus  T_\smooth$ as before.

But in order to make a connection with the previous picture and with
the geometry of toric hypersurfaces we need to twist the complex
structure on $W_\epsilon(\lambda)$ by the element of
$H^1(\smooth,\T^{d-1})$ associated with the phase parameters
$\theta_v$.
The resulting K\"ahler manifold $W_\epsilon(\lambda,\theta)$ can be
canonically identified with $W(\epsilon^{-1}\lambda,\theta)$
constructed  by the first method starting with the $\epsilon$-rescaled
K\"ahler affine structure.

\section{Vector fields, foliations and K\"ahler potentials}

This section describes two important ingredients for a later
consideration of the geometry of toric hypersurfaces. The first is the
foliations of the amoebas which will induce the torus fibrations of
the hypersurfaces. The second is a K\"ahler potential on the toric
variety which induces the metric.

From now on we will fix a projective family of bi-PIKAS on
$\smooth$. We will refer to a member of type $(\lambda,\nu)$ as the
$(\lambda,\nu)$-bi-PIKAS.
Also we will fix some linear functional $\ell$ positive on the
secondary cone $\SC(S)$. This will allow us to talk about the scale
$\ell(\lambda)$ of the vector $\lambda$. We will abuse the notation
and denote the analogous scale for $\nu$ by $\ell(\nu)$.

Let denote by $\overrightarrow{m}$ the vector in $\R^{\dD_\Z}$ with
$\overrightarrow{m}(m)=1$ and $0$ otherwise, and let
$\overrightarrow{\mathbbm{1}}:=\sum_{m\in\dD_\Z}\overrightarrow{m}$.
For a real number $\beta$ we will often write $\lambda+\beta$ meaning
$\lambda+\beta\cdot\overrightarrow{\mathbbm{1}}$.
We will say that $\beta>0$ is {\it small in the $\lambda$-scale}, or
simply $\lambda$-small, if the vectors
$\lambda+\beta(\sum\pm\overrightarrow{m_i})$, for all possible
collections $\{m_i\}\subset\dD_\Z$, are still in the interior of the
secondary cone $\SC(S)$. And similar for the $\nu$-scale.

\subsection{Neighborhoods of the discriminant}
For a given $(\lambda,\nu)$ we define subsets $U_v^\beta\subset U_v$
and $V_w^{\beta^\vee}\subset V_w$ whose union will give the complement
of a neighborhood of $D$ depending on two real parameters
$\beta,{\beta^\vee} >0$. We assume $\beta,{\beta^\vee}$ to be small in
the $\lambda,\nu$ scales, respectively. In what follows we identify
$\Sigma$ with $\dDS$ and $\dDT$ via the maps $\phi$ and $\hat\phi$
associated with the bi-PIKAS of type $(\lambda,\nu)$.

For $v\in\vertS$ we let $U_v^\beta$ be the set of points in the
corresponding facet of $\DSv$ which lie in the closed polyhedron
$Q^\lambda_{(v\mid\{0\})}(\beta)$ (cf. \cite[Section~3.2]{HZh}), and
similar for $V_w^{\beta^\vee}$. Explicitly,
\begin{gather*}
  U_v^\beta:=\{n\in U_v \suchthat \<m,n\>+\lambda(m)\le -\beta, \text{
    all } m\in\D_\Z\setminus \{v,0\}\},\\
  V_w^{\beta^\vee}:=\{m\in V_w \suchthat \<m,n\>+\nu(n)\le
  -{\beta^\vee}, \text{ all } n\in\Dv_\Z\setminus \{w,0\}\}.
\end{gather*} 
Because $\beta,\beta^\vee$ are small the sets
$U_v^\beta$ and $V_w^{\beta^\vee}$ are non-empty.  
We define the smooth part of $\Sigma$ as
$$\Sigma^\sm:=\bigcup_{v\in \vertS} U_v^\beta \cup 
  \bigcup_{w\in \vertT} V_w^{\beta^\vee}.
  $$
Then the neighborhood of $D$ is defined as  the complement to all
these closed sets in $\Sigma$:
\begin{equation*}
  N_{\lambda,\nu}^{\beta,\beta^\vee}(D):=\Sigma\setminus \Sigma^\sm. 
\end{equation*}
An important observation is that
$N_{\lambda,\nu}^{\beta,\beta^\vee}(D)\to D$ as $\beta,\beta^\vee\to
0$ in the scales of $\lambda,\nu$, respectively.

\subsection{Regularization}\label{sec:regularization}
We will be smoothing various functions later on in this section so let
us recall the standard regularization techniques. Given a convex
bounded domain $P\subset\R^n$ with piece-wise smooth boundary let
$\rho$ be a mollifier whose support is $P$, i.e. a positive
$\C^\infty(\R^n)$ function vanishing exactly outside $P$ and such that
$\int_{\R^n}\rho dx =1$. For instance, if $P$ is a polytope in $\R^n$
given by a collection of inequalities $\{\<v_i,x\>+\lambda_i\leq 0\}$
one can take the usual bell-shaped function
\begin{equation*}
  \rho=c\prod_i \rho_i, \quad \text{ where }
  \rho_i(x)=
  \begin{cases}
    e^{\frac1{\<v_i,x\>+\lambda_i}},& \quad  \<v_i,x\>+\lambda_i\leq 0\\
    0,&\quad \<v_i,x\>+\lambda_i\geq 0
  \end{cases},
\end{equation*}
and the constant $c$ is determined from the normalization.
Set $\rho_{h}:=\frac1{h^n}\rho(\frac{x}{h})$. For any locally
integrable function $u\in L^1_{loc} (\R^n)$ we can apply the standard
regularization procedure by taking the convolution with $\rho_h$:
$$ u_h(x) := \rho_h \ast u =\int_{\R^n} \rho_h(x-y) u(y)\, dy.$$
Then, if $u\in C^p(\R^n)$, the functions $u_h$ are $\C^\infty$ and
approaching $u$ as $h\to 0$ in the $C^p$-norm uniformly on any compact
in $\R^n$.

Below we list some elementary properties of $u_h$ which will be useful
later.
\begin{prop}\label{prop:convexity}
  Let $\Omega$ be a convex domain in $\R^n$, then
  \begin{itemize}
  \item $u_h$ is linear in a $v$-direction in $\Omega$ if $u$ is
    linear in the $v$-direction in the Minkowski sum $\Omega+h(-P)$,
    with $\<v,\nabla u_h\>=\<v,\nabla u\>$.
  \item $u_h$ is convex in $\R^n$ if $u$ is. The gradient $\nabla
    u_h(x), x\in\Omega$, is always inside the convex hull of all
    possible gradients of $u$ in $\Omega+h(-P)$.
  \item $u_h$ is strictly convex in $\Omega+h(-P)$ if $u$ is convex in
    $\R^n$ and strictly convex in $\Omega$.
  \end{itemize}
\end{prop}
\begin{proof}
  All statements are simple consequences of the following
  observation. The convolution $\rho_h\ast u$ is a weighted averaging
  of $u$ over the (translated) support of $\rho$. The same holds for
  all derivatives of $u$ as well.
\end{proof}

One can also apply the regularization procedure by taking the
convolution with the mollifier parameter depending (smoothly) on the
point $x\in\R^d$. That is $u_{h}(x) :=\int_{\R^n} \rho_{h(x)}(x-y)
u(y)\, dy$. Or, even, more generally the entire shape of the support
of the mollifier can smoothly depend on the center of convolution. We
have used this technique of varying support to prove the existence of
bi-PIKAS in Proposition~\ref{prop:existence}.

\subsection{Fibration}
We will construct a foliation of  $\R^d \setminus \DSv$ by straight
lines/rays by specifying a ``convex'' vector field on $\dDS$, which is
smooth in $\Sigma^\sm$.

Recall that the piece-wise linear functions $L_\nu:(\R^d)^*\to\R$ and
$L_\lambda:\R^d\to\R$, the Legendre transforms of $\nu$ and $\lambda$,
were defined as
$$L_\nu(m):=\max_{n\in\Dv_\Z}\{\<m,n\>+\nu(n)\}, \quad
L_\lambda(n):=\max_{m\in\D_\Z}\{\<m,n\>+\lambda(m)\}.
$$
We fix a mollifier $\rho^\vee$ on $(\R^d)^*$ with support in $\D$ and
consider smooth functions $L_{\nu,h^\vee}:=\rho^\vee_{h^\vee}\ast
L_{\nu-h^\vee}$. From the properties of regularization (for ${h^\vee}$
small in the $\nu$-scale) the slopes $\nabla L_{\nu,{h^\vee}}(x)$
always lie in $\dDv$, for any $x$ not in the interior of $\DT$. In
particular, the gradient of $L_{\nu,{h^\vee}}$ gives a map $\nabla
L_{\nu,{h^\vee}}:\dDT\to\dDv$.

Now we can use the identification of $\dDT$ with $\dDS$ via
$\phi\hat\phi^{-1}$ given by the $(\lambda,\nu)$-bi-PIKAS to define a
$\dDv$-valued vector field $\X_{h^\vee}$ on $\dDS$ by
$$\X_{h^\vee}(q):=\nabla L_{\nu,h^\vee}(\phi\hat\phi^{-1}q).$$ 
We can extend this vector field $\X_{h^\vee}$ to $\R^d \setminus \DSv$
using the following identification of $\dDT$ with $\dDv_{\lambda+t}$,
for any $t\geq 0$. If $\hat\Phi$ is the (dual)
$(\lambda,\nu)$-bi-PIKAS potential, then $\hat\Phi+tL_{\nu,{h^\vee}}$
is a strictly convex function on $\DT$. In particular, its gradient
defines a bijection $\nabla\hat\Phi+t\nabla L_{\nu,{h^\vee}}:
\dDT\to\dDv_{\lambda+t}$. We can write the vector field $\X_{h^\vee}$
on $\R^d \setminus \DSv$ explicitly, using the fact that
$$L_\lambda(x)=t \Leftrightarrow x\in\dDv_{\lambda+t}, \text{ for any
} t\geq 0.$$
Namely,
$$\X_{h^\vee}(x):=\nabla
L_{\nu,h^\vee}([\nabla\hat\Phi+L_\lambda(x)\nabla
L_{\nu,{h^\vee}}]^{-1}(x)).
$$
We summarize properties of this vector field in the following lemma: 
\begin{lemma}\label{lemma:fibration}
For any $h^\vee>0$, small in the $\nu$-scale, the vector field
$\X_{h^\vee}$ induces a (straight line/ray) foliation $\F_{h^\vee}$ of
$\R^d \setminus \DSv$, smooth over $\smooth$, such that
\begin{enumerate}
\item $\X_{h^\vee}(q)=w$ for $q\in V_w^{\beta^\vee}$.
\item The value of $\X_{h^\vee}(q)$ is in
  $(\carrier\sigma)^\vee\subset\dDv$ for $q \in
  F_\sigma\subset\dDS$. In particular, $\<v,\X_{h^\vee}\>=1$ in
  $U_v$.
\item For any $q\in\smooth$,
  $\left|\nabla\X_{h^\vee}(q)\right|\le\frac {C|g_{ij}(q)|}{h^\vee}$,
  where $C$ is a constant independent of $\lambda$ and $\nu$, and the
  gradient and the metric $g_{ij}$ are taken in affine coordinates.
\end{enumerate}
\end{lemma}
\begin{proof}
  From the definition it is easy to see that
  $\X_{h^\vee}(q+t\X_{h^\vee}(q))=\X_{h^\vee}(q)$, for $q\in\dDS$, which
  immediately implies the straight ray foliation.
  
  Note that as ${h^\vee}\to 0$, $L_{\nu,{h^\vee}}$ converges to $L_\nu$
  uniformly in $(\R^d)^*$. For $h^\vee=0$ the (discontinuous) vector
  field $\X$ has (discrete) values in $\vertT$, and (1) and (2) follow
  immediately from the combinatorics of $\Sigma$. They remain true after
  regularization as well, which is guaranteed by the
  Proposition~\ref{prop:convexity}.
  
  The bound (3) on the derivatives of $\X$ follows from a standard
  estimate for regularization of piece-wise smooth function $L_\nu$. In
  the {\em dual} affine coordinates $\Hess L_\nu$ is a Dirac
  $\delta$-like distribution supported on $\dV$. The norm of its
  convolution with $\rho$ is bounded by $\frac C{(h^\vee)^k}$, where $k$
  is the codimension of the support. The constant $C$ takes into account
  the combinatorics of the polytope $\D$, the particular form of the
  mollifier $\rho$ and the choice of the norm on $\R^{d-1}\cong
  T_q(\smooth)$. The metric $g_{ij}$ appears from the chain rule:
  $\nabla\X_{h^\vee}=g\cdot \Hess L_{\nu,h^\vee}$.
  
  Finally, the smoothness of the vector field $\X_{h^\vee}$, and hence
  the smoothness of the foliation $\F$, follow from smoothness of the
  map $\phi\hat\phi^{-1}$ on $\smooth$.
\end{proof}

\subsection{K\"ahler metrics on the toric variety}
First, we would like to extend the bi-PIKAS potential to $\R^d$ by taking the Legendre transform of $\hat\Phi$. Namely, 
$$\Phi(x)=\max_{y\in\DT}\{\<x,y\>-\hat\Phi(y)\}.$$ 
Similarly, we extend $\hat\Phi$ to a function on $(\R^d)^*$. We will abuse the notation $\Phi,\hat\Phi$ for the extended potentials.

$\Phi$ is a $C^1$-function, smooth when restricted to any strata of $\DSv$. Its Hessian $\Hess\Phi$ is continuous at $\dU$, but blows off at $\dV$. And something drastic happens at the discriminant $D=\dU\cap\dV$. 

Next we regularize the $C^1$-potential $\Phi$ to get a smooth convex
function on $\R^d$ which we will use  later on to define a K\"ahler
potential on the toric variety $X_T$. Let $\rho$ be a mollifier with
support in $-\Dv$. We define $\Phi^\sm_h:=\rho_h\ast\Phi.$

\begin{remark}
  The constructed vector field, foliation, potential, etc., depend on
  the pair $(\lambda,\nu)$, as well as on the regularization
  parameters $h,h^\vee$. But to simplify the notations for
  $\X,\F,\Phi,\Phi^\sm$ we will often leave only those indices which
  are important in a current consideration and omit the rest when
  there is no confusion possible.
\end{remark}

Before constructing a K\"ahler potential on the toric variety $X_T$ we
need another technical statement.
\begin{lemma}\label{lemma:potential}
  For $\tau \in T$, the $\tau$-slope of $\Phi^\sm_{h}$ is equal to the
  $\tau$-slope of $-\nu|_\tau$ in some translation $R_{\geq\tau}$ of
  $\cup_{\tau'\geq\tau}\cone(\tau')$.
\end{lemma}
\begin{proof}
  Consider $h=0$ first. For a simplex $\tau\in T$, let $F_\tau$ be the
  corresponding face of $\DT$, and we set $R_\tau:=
  \phi\hat\phi^{-1}(F_\tau)+\cone(\tau)$. Then the set $ \cup_{\tau'
    \ge \tau} R_{\tau'}$ contains some translation $R_{\ge \tau}$ of
  $\cup_{\tau'\geq\tau}\cone(\tau')$.

  On the other hand, $\Phi$ is the Legendre transform of
  $\hat\Phi|_\DT$.
  Hence, if $\nabla\Phi(x)=y$ and $r$ is in the normal cone to $\DT$ at
  $y\in\DT$, then $\nabla\Phi(x+r)=\nabla\Phi(x)=y$.
  So we see that for $x\in R_\tau$ the gradient $\nabla\Phi(x)$ takes
  values in the face $F_\tau$ of $\DT$ because $x$. In particular, the
  $\tau$-slopes of $\Phi$ are equal to $-\nu|_\tau$.

  For $h>0$ the statement of the lemma follows from the case $h=0$ and
  the Proposition~\ref{prop:convexity}. The translated cones
  $R_{\ge\tau}$ become shifted into their interiors by some vectors of
  size $h$.
\end{proof}

Now we can use $\Phi^\sm=\Phi^\sm_h$ with any $h>0$ to define a
K\"ahler potential on $X_T$. For an element $z=\{z_1,\dots,z_d\}\in\Cstar^d$ we will use the notations 
$$\log|z|:=\{\log|z_1|,\dots,\log|z_d|\}\in \R^d,\quad
\Arg(z):=\frac1{2\pi}\{\arg(z_1),\dots,\arg(z_d)\}\in\T.$$

\begin{prop}
  The $(1,1)$-form defined on $(\Cstar)^d$ by
  $$\eta:=\frac{\sqrt{-1}}{2\pi}\partial\bar{\partial}
  \Phi^\sm\left(\log|z|\right)$$
  extends to a smooth (in the orbifold sense) non-negative definite
  $(1,1)$-form on $X_T$ in the cohomology class $[\eta]=[\nu]$.
\end{prop}

\begin{proof}
  First, we rewrite the form $\eta$ on $(\Cstar)^d$ as
  \begin{equation*}
    \eta=\frac{\sqrt{-1}}{2\pi}\partial\bar{\partial}   
    \Phi^\sm\left(\log|z|\right)
    =\left\<  
      d\left(\nabla\Phi^\sm\left(\log|z|\right)\right) \wedge d\Arg(z)
    \right\>,
  \end{equation*}
  where "$\<\ \wedge\ \>$" means also the $\<\ ,\ \>$-pairing between
  the $(\R^d)^*$-valued gradient $\nabla\Phi^\sm$ and the
  $\R^d$-valued 1-form $d\Arg(z)$.

  For a simplex $\tau\in T$ we want to show that $\eta$ extends to the
  toric subvariety $Z_\tau$ associated to $\tau$. If $X_T$ is smooth,
  then in a neighborhood of $Z_\tau$ we can choose the coordinates
  similar to those from \cite[Lemma~3.9]{HZh}. That is, we choose a
  basis $\{e_i\}$ such that
  \begin{equation*}
    \<e_i, w_j\>=-\delta_{ij}, i=1,\dots,\dim\tau+1 \text{ and }
    \<e_i,\tau\>=0, i=\dim\tau+2,\dots,d.
  \end{equation*} 
  Then, in the coordinates $y_i=z^{e_i}$ the equations for the
  subvariety $Z_\tau\subset X_T$ are $y_i=0, i=1,\dots,\dim\tau+1$.

  According to the theory of toric varieties (cf., e.g. \cite{Fulton})
  a neighborhood of the toric subvariety $Z_\tau$ lies in the closure
  of $\log^{-1}(R_{\ge\tau})$, where $R_{\ge\tau}$ is any translation
  of the cone $\cup_{\tau'\geq\tau}\cone(\tau')$. But by the
  Lemma~\ref{lemma:potential} the directional derivatives
  $\<\nabla\Phi^\sm,w_i\>$, $w_i\in\tau$, are constant in some
  translation $R_{\ge\tau}$ of  $\cup_{\tau'\geq\tau}\cone(\tau')$. Hence, in a
  neighborhood of $Z_\tau$ the form $\eta_a$ written in the above
  coordinates is independent of $y_i, i=1,\dots,\dim\tau+1$, and,
  thus, can be extended to $Z_\tau$.

  In case when $X_T$ is an orbifold we may not be able to choose an
  integral basis $\{e_i\}$ with the above conditions. This corresponds
  to the fact that we may need to go to a finite cover to get a smooth
  form by weakening the first set of conditions to be
  $\<e_i,\tau\>\in\Z$. But the rest of the argument goes through. 

  Finally, the cohomology class of a $\T$-invariant $(1,1)$-form on a
  complete toric variety is determined by the image of its moment
  map. But the moment map for $\eta$ is given on $(\Cstar)^d$ by
  \begin{equation*}
    \mu(z)=\nabla\Phi^\sm\left(\log|z|\right),
  \end{equation*}
  whose extension to the whole toric variety $X_T$ has the image
  $\DT$. Hence the class of $\eta$ is $[\nu]$.
\end{proof}

Finally we can add to $\eta$ a (small) positive multiple of a
K\"ahler (e.g., Fubini-Study) form $\omega_0$. Thus we get a true
K\"ahler form $\omega=\eta+\epsilon\omega_0$ on $X_T$ in the class
$[\nu] +\epsilon[\omega_0])$.

\section{Geometry of Calabi-Yau toric hypersurfaces} 

A Calabi-Yau hypersurface $Z_a$ is given by the closure of the set
\begin{equation*}
  Z_a^\aff :=\{z \in (\C \setminus \{0\})^d \suchthat
  \sum_{m\in\D_\Z\setminus \{0\}} a_m z^m = 1 \}
\end{equation*}
in the toric variety $X_T$. 
From now on we set $\lambda:=\log|a|$ and require it to be in a proper
subcone of the secondary cone $\SC(S)$. Also, for
non-zero $a_v$, we set $\theta_v:= \frac1{2\pi}\arg (a_v)$.

\subsection{An embedding of $\Hsm$ into the model torus bundle}

In \cite{HZh} we have used the GKZ machinery \cite{GKZ} for the
monomial estimates in the equation of $Z_a$ to establish an embedding
of $\Hsm$ into $W_a:=W(\log|a|,\Arg(a))$. The same estimates  can be
used to find bounds on the discrepancy of this embedding from being
holomorphic and isometric. Establishing these bounds will occupy the
rest of the section.

To define the fibration we assume that $\lambda=\log|a|$ is
sufficiently far in the interior of $\SC(S)$, so that
$c:=\log|\D_\Z|$ is small in the $\lambda$-scale.
\begin{lemma}\label{lemma:amoeba}
  The amoeba $\am=\log(Z_a^\aff)$ lies outside of
  $\Dv_{\lambda-c}$. In particular, the foliation
  $\F_{\lambda-c,\nu,h}$ of $\R^d\setminus\Dv_{\lambda-c}$  induces
  the fibration $\am\to\dDv_{\lambda-c}$ by projection along the
  leaves.
\end{lemma}
\begin{proof}
  Note that for any $m\in\D_\Z\setminus\{0\}$, if
  $\<m,\log|z|\>+\log|a|\leq -\log|\D_\Z|$, then $|a_m z^m|\leq
  \frac 1{|\D_\Z|}$. Hence the equation
  $\sum_{m\in\D_\Z\setminus\{0\}}a_mz^m=1$ cannot have solutions for
  $x\in\log^{-1}( \Dv_{\lambda-c})$.
\end{proof}
From now on we will use the $(\lambda-c,\nu)$ bi-PIKAS to identify
$\Sigma$ with $\dDv_{\lambda-c}$ and fix the vector field and the
foliation in $\R^d\setminus\Dv_{\lambda-c}$. With this identification,
given a subset $U\in\Sigma$ we denote by $X(U)$ the closure of the set
$\log^{-1}(\cup_{q\in U} \F_q)$ in the toric variety $X_T$
(cf. \cite{HZh}).
Then the smooth part $\Hsm$ of the hypersurface is defined as:
$$ \Hsm:=Z_a\cap X(\Sigma\setminus
N_{\lambda-c,\nu}^{\beta,\beta^\vee}(D)).$$

We define the map $\psi$ over the charts $V_w^{\beta^\vee}$ as the
restriction to $Z_a$ of the quotient map:
\begin{equation*}\label{eq:projection}
  z\mapsto (\log|z|/w, \Arg(z)/w) \in (\R^d/w,\T/w). 
\end{equation*}
Note that if $x$ is in a boundary toric divisor $Z_w, w\in \vertT$,
then $\log|z|$ and $\Arg(z)$ are not well defined, but $\log|z|/w$ and
$\Arg(z)/w$ are. Hence, the map $\psi$ is well defined over
$X(V_w^{\beta^\vee})$. The meaning of this map is the choice of local
coordinates for $Z_a$ near $z$ (cf. \cite[Lemma 3.9]{HZh}). Hence it
is holomorphic.

Before defining the map over the charts $U_v^\beta$ let us first make
some estimates in the spirit of Lemma~\ref{lemma:amoeba}. For an
element $\theta\in\R/\Z$ we will write $|\theta|<c$ if the (standard
Euclidean) distance from $\theta$ to 0 is less than $c$. This
inequality is vacuous for $c\ge 1/2$.
 
\begin{lemma}\label{lemma:nbhd}
If a point $x$ of $Z_a$ lies in $X(U_v^\beta)$, then
\begin{gather*}  
  |\<v,\log|z|\>+\lambda(v)|\le C(\beta)e^{-\beta}\\
  |\<v,\Arg(z)\>+\theta_v|\le C(\beta)e^{-\beta},
\end{gather*}
where $C(\beta)\to \log|\D_\Z|$  as $\beta\to \infty$.  
\end{lemma}
\begin{proof}
  First of all note that since $X(U_v^{\beta})\cap Z_a
  \subset (\Cstar)^d$, both $\log|z|$ and $\Arg(z)$ are well
  defined. Also by (2) of Lemma~\ref{lemma:fibration} the values of
  the vector field $\X$ are in $(\carrier v)^\vee$. Hence, for any
  ${q\in U_v^{\beta}}$ the ray $\F_q$ is in $
  Q^\lambda_{(v|\{0\})}(\beta)$ and the standard estimates on values
  of the monomials at $z\in\log^{-1}(Q^\lambda_{(v|\{0\})}(\beta))$
  apply:
  $$|a_m z^m|\leq e^{-\beta}|a_v z^v|, \text{ all } m\neq v,\{0\}.$$
  Or putting them all together we have
  \begin{gather*}
    \left|\frac1{a_vz^v }\sum_{m\ne v,\{0\}}a_mz^m\right|\leq
    |\D_\Z|e^{-\beta}.
  \end{gather*}
  Hence,
  \begin{gather*}
    \left| \log\left( 1+\frac1{a_vz^v }\sum_{m\ne
          v,\{0\}}a_mz^m\right)\right|\leq
    C'(\beta)\cdot|\D_\Z|e^{-\beta},
  \end{gather*}
  where $C'(\beta)\to 1$ as $\beta\to\infty$. Writing the equation of
  $Z_a$ in $X(U_v^\beta)$ as
  \begin{equation*}
    a_vz^v \left(1+\frac1{a_vz^v }\sum_{m\ne v,\{0\}}a_mz^m\right) =1
  \end{equation*}
  or, equivalently,
  \begin{equation*}
    \log (a_vz^v) =-\log\left( 1+\frac1{a_vz^v }\sum_{m\ne
        v,\{0\}}a_mz^m\right)
  \end{equation*}
  will give the claimed estimates.
\end{proof}

Now, identifying the tangent spaces of the torus fibers
$\T_x=\log^{-1}(x)$ with $\R^d$, we can pull back the vector field
$\X(x)$ to get a (constant) vector field on $\T_x$. Then the map
$\psi$ for the points in $X(U_v^\beta)\cap Z_a \subset (\Cstar)^d$
will be defined as:
\begin{equation*}
  \begin{split}
    \psi(z):=
    &\bigl(\log|z|-(\<v,\log|z|\>+\lambda(v))\X(\log|z|),\\ 
    &\Arg(z)-(\<v,\Arg(z)\>+\theta_v)\X(\log|z|)\bigr) \in
    \bigl(\R^d_v(\lambda(v), \T_v(\theta_v)\bigr).
  \end{split}
\end{equation*}  
Here if $\beta$ is large enough, i.e. $C(\beta)e^{-\beta}<\frac12$,
then according to the Lemma~\ref{lemma:nbhd}  there is a preferred
{\em continuous} lift of  $\<v,\Arg(z)\>+\theta_v$ to (the
neighborhood of 0 in) $\R$. We use this lift to first define the value
for $(\<v,\Arg(z)\>+\theta_v)\X$ in $\R^d$ and then project it back to
$\T=\R^d/\Z^d$.

Since $\X_q=w$ for $q\in V_w^{\beta^\vee}$ the definition of the map
$\psi$ over the charts $V_w^{\beta^\vee}$ is consistent with the above
definition on possible overlaps $U_v^\beta\cap
V_w^{\beta^\vee}$. Thus, we have a well defined map $\psi:\Hsm\to W_a$
which is an embedding \cite{HZh}.

\subsection{Estimates on complex structures and metrics}

Let $J_{Z_a}$ and $J_{W_a}$ denote the complex structure operators on
the tangent spaces to $Z_a$ and $W_a$ respectively. We would like to
say that the embedding $\psi$ is holomorphic up to a small order
terms.

We have already mentioned that $\psi$ is {\em precisely} holomorphic
over the charts $V_w^{\beta^\vee}$. To measure the discrepancy at
$x\in U_v^\beta$ we will fix some (Euclidean) norm on $\R^{d-1}$ (they
are all equivalent) to induce a norm on the tangent space
$\R^{d-1}\oplus\R^{d-1}\cong T_x{W_a}$. Let $C_{\lambda,\nu}(\beta)$
be the (uniform on $U_v^\beta$) bound for the $(\lambda,\nu)$ bi-PIKAS
metric $g_{\lambda,\nu}$ written in the affine coordinates in
$U_v^\beta$.

\begin{lemma}\label{lemma:complex}
  As $\beta\to\infty, \beta^\vee\to 0$, the linear map 
  $$d\psi\circ J_{Z_a}\circ(d\psi)^{-1} -J_{W_a}:T_x{W_a}\to
  T_x{W_a}$$
  is of order $C_{\lambda-c,\nu}(\beta)\cdot
  O(\frac1{\beta^\vee}e^{-\beta})$.
\end{lemma}
\begin{proof}
  First we apply estimates similar to those in Lemma~\ref{lemma:nbhd}
  to the differential $d\psi$, which we think of as an element in
  $(\R^d\oplus\R^d)\otimes((\R^d)^*\oplus(\R^d)^*)$.
  \begin{equation*}
    d\psi=
    \left(\begin{array}{cc}
        \mathbbm{1}-\X\otimes v-(\<v,\log|z|\>+\lambda(v))\nabla \X &
        0 \\
        -(\<v,\Arg(z)\>+\theta_v)\nabla \X & \mathbbm{1}-\X \otimes v 
      \end{array}\right)
    \left(\begin{array}{c} d\log|z| \\ d\Arg(z) \end{array}\right)
  \end{equation*}
  Note that the complex structures $J_{Z_a}$ and $J_{W_a}$ would match
  exactly via $d\psi$ if there were no $\nabla\X$ terms (this is what
  happens in the charts $V_w^{\beta^\vee}$ where $\X$ is constant).
  
  According to (2) of Lemma~\ref{lemma:fibration} $\<v,\X\>=1$ in
  $U_v^\beta$. Hence the projection operator
  \begin{equation*}
    \left(\begin{array}{cc}
        \mathbbm{1}-\X\otimes v & 0 \\
        0 & \mathbbm{1}-\X \otimes v 
      \end{array}\right)
  \end{equation*}
  has a norm of order 1 when restricted to $Z_a$, and the desired
  bound on $d\psi\circ J_{Z_a}\circ(d\psi)^{-1} -J_{W_a}$ will follow
  from estimating the $\nabla\X$ terms.

  But according to the Lemma~\ref{lemma:nbhd} we have the uniform
  bounds:
  $$
  |\<v,\log|z|\>+\lambda(v)|\le Ce^{-\beta},\
  |\<v,\Arg(z)\>+\theta_v|\le Ce^{-\beta}.
  $$
  On the other hand, by (3) of Lemma~\ref{lemma:fibration} the
  gradient $\nabla \X(\log|z|)$ is bounded by
  $O(\frac1{\beta^\vee})\cdot
  \left|g_{\lambda-c,\nu}(\log|z|)\right|$.
\end{proof}

The K\"ahler form on $\Hsm$ is defined as the restriction of the
K\"ahler form $\omega$ on $X_T$ which in $(\Cstar)^d$ is given by:
$$\omega=\frac{\sqrt{-1}}{2\pi}\partial\bar{\partial}
\Phi^\sm_{\lambda+\gamma,\nu,h}\left(\log|z|\right)+
\epsilon\omega_{0},
$$
where $\omega_{0}$ is a fixed (e.g., the Fubini-Study) K\"ahler
form.
We will compare the metric induced by $\omega$ with the (degenerate)
scalar product on $W_a$ induced by the $(\lambda,\nu)$-bi-PIKAS.

To make these estimates we will need to introduce some bounds (in a
Euclidean metric in $\R^d$) all of which follow essentially from  the
definition of the bi-PIKAS family:
\begin{equation*}
  \begin{split}
    &\left|g_{ij}(x)\right|<C_0(\beta), \text{ from
      Lemma~\ref{lemma:complex} above, uniformly for } x\in
    U_v^\beta,\\
    &\left|\Hess \Phi_{\lambda+\gamma}(x)\right|<C_1(\gamma), \text{
      uniformly for } x\in\Dv_{\lambda+\gamma/2},\\
    &\left|g_{\lambda+\gamma,\nu}(q)-g_{\lambda,\nu}(q)
    \right|<C_2(\gamma), \text{ uniformly for } q\in\dU,\\
    &\left|\Hess \Phi^\sm_{\lambda+\gamma,h}(x)-\Hess
      \Phi_{\lambda+\gamma}(x)\right|<C_3(h), \text{ uniformly for }
    x\in\Dv_{\lambda+\gamma/2},\\
    &\left|\Hess \Phi_{\lambda+\gamma}|_{\R^d_v}(q+t\X(x))-\Hess
      \Phi_{\lambda}|_{\R^d_v}(q)\right|<C_4(\gamma), \text{ uniformly
      for }\\
    & \hskip 1in  -\gamma<t<\gamma, q\in U_v^\beta\subset\dDS,\\
    &\left|\Hess \Phi_{\lambda}(x)-\Hess K_w(q)\right|<C_5(\beta,c),
    \text{ uniformly for }\\
    & \hskip 1in q\in N^\beta(\dU)\cap V_w^{\beta^\vee}, x\in
    \F_q:=q+\R_{\geq-c}\cdot w,
  \end{split}
\end{equation*} 
where $C_0(\beta)\to \infty,\ C_1(\gamma)\to\infty,\ C_2(\gamma)\to
0,\ C_3(h)\to 0,\ C_4(t)\to 0,\ C_5(\beta,c)\to 0$ as (all of) the
corresponding parameters go to 0. The last inequality follows from
$\Phi_{\lambda}\in C^2(\dU\setminus D)$, where the local potential
$K_w(x)$ is pulled back from the quotient.

\begin{lemma}\label{lemma:scalar}
  Under the embedding $\psi:\Hsm\to W_a$ the scalar products agree up
  to terms of order $C_2(\gamma) + C_3(h) + C_1(\gamma)
  \frac{C_0(\beta)}{\beta^\vee} e^{-\beta} + C_4(\gamma) +
  C_5(\beta,c) + O(\epsilon).$
\end{lemma}
\begin{proof}
  Let $x\in X(U_v^\beta)\cap Z_a$. Assuming $e^{-\beta}<\gamma/2$ we
  have
  $$\left|\Hess \Phi^\sm_{\lambda+\gamma,h}(x)-\Hess
    \Phi_{\lambda+\gamma}(x)\right|<C_3(h).$$
  The difference between $\Hess \Phi_{\lambda+\gamma}|_{T_xZ_a}$ and
  $\Hess \Phi_{\lambda+\gamma}|_{T_{\psi(x)}W_a}$ consists of two
  terms. The first term $C_4(\gamma)$ appears from comparing $\Hess
  \Phi_{\lambda+\gamma}|_{\R^d_v}$ at $q+t\X(x)$ with $\Hess
  \Phi_{\lambda}|_{\R^d_v}$ at $q$. The other term
  $C_1(\gamma)\frac{C_0(\beta)}{\beta^\vee}e^{-\beta}$ reflects the
  error in the alignment of the tangent spaces via the map $d\psi$ in
  the proof of Lemma~\ref{lemma:complex}.

  Now let $x\in X(V_w^{\beta^\vee})$. We will just need to check
  points in $\F_q$ for $q\in V_w^{\beta^\vee}\setminus \bigcup
  U_v^\beta= N^\beta(\dU)\cap V_w^{\beta^\vee}$. Note that $\Hess
  \Phi_{\lambda}$ is bounded in $\F(V_w^{\beta^\vee})$ and is
  continuous at $\dU\setminus D$ (in particular, the Hessian vanishes
  into the $w$-direction at $\dU\setminus D$).  Hence the $C_5$ bound
  from above are also valid for the regularization $\Hess
  \Phi^\sm_{\lambda+\gamma,h}$. Namely,
  $$ \left|\Hess \Phi^\sm_{\lambda+\gamma,h}(x) -
    g_{\lambda+\gamma,\nu}(q) \right| < C_5(\beta,c),$$
  with possibly different function $C_5$. The discrepancy between
  $g_{\lambda+\gamma,\nu}$ and $g_{\lambda,\nu}$ in $N^\beta(\dU)\cap
  V_w^{\beta^\vee}$ is encoded in the $C_2(\gamma)$ term.

  Finally, the term $\epsilon\omega_{0,\nu}$ in $\omega$ can be
  bounded by $O(\epsilon)$.
\end{proof}

\subsection{The Gromov-Hausdorff limits of one-parameter families}
We will apply the results of the previous sections to the situation
considered in \cite{HZh} to draw a consequence mostly related to the
mirror symmetry conjecture. Let $\lambda_0$ be an integral vector in
the interior of the secondary cone $\SC(S)$. We consider an
1-parameter family of the hypersurfaces $Z_s$ defined as closures in
$X_T$ of
\begin{equation*}
 Z_s^\aff :=\{z \in (\C \setminus \{0\})^d \suchthat
 \sum_{m\in\D_\Z\setminus \{0\}} a_m s^{\lambda_0(m)} z^m = 1\}.
\end{equation*}
Choose an integral vector $\nu_0$ in the interior of $\SC(T)$ and
consider $\Sigma$ with the metric space structure given by the
bi-PIKAS $(\lambda_0,\nu_0)$.

Also consider an one-parameter family of (non-compact) K\"ahler
manifolds $W_s:=W(\log|a|+\log|s|\cdot\lambda_0,
\Arg(a)+\Arg(s)\cdot\lambda_0)$, whose metric and complex structure
are induced from the
$\left(\log|a|+\log|s|\cdot\lambda_0,\frac{\nu_0}{\log|s|}\right)$
bi-PIKAS.

\begin{thm} 
  As $|s|\to\infty$ one can choose smooth portions of the
  hypersurfaces $Z_s^\sm\subset Z_s$, the embeddings $\psi_s:
  Z_s^\sm\hookrightarrow W_{s}$ and a family of K\"ahler metrics on
  $Z_s$ in the class $\frac{\nu_0}{\log|s|}(1+o(1))$ such that the
  pairs $(Z_s,Z_s\setminus Z_s^\sm)$ converges to the pair
  $(\Sigma,D)$ in the Gromov-Hausdorff sense, and the maps $\psi_s$
  identify (uniformly in $x\in Z_s^\sm$) the scalar products and the
  complex structures on the tangent spaces $T_x Z_s$ and
  $T_{\psi_s(x)} W_{s}$ up to terms of order $o(1)$.
\end{thm}
\begin{proof} 
  We consider the bi-PIKAS family in a neighborhood of
  $(\lambda_0,\nu_0)$ and extend it by rescaling to
  (a neighborhood of) the ray $(\lambda,\nu) = \left( \log|a| +
    \log|s| \cdot \lambda_0, \frac{\nu_0}{\log|s|} \right)$  
  in $\SC(S)\times\SC(T)$.

  The $C_i$-estimates for the $(\lambda,\nu)$ bi-PIKAS considered in
  the bi-PIKAS metric (rather than in Euclidean) are equivalent to the
  corresponding estimates for the rescaled structure
  $((\log|s|)^{-1}\lambda,\log|s|\nu)$ made in the Euclidean metric as
  before. This is because the Euclidean metric on $\dD_{\lambda_0}$ is
  equivalent to the bi-PIKAS metric on $\Sigma$. But to pass from
  $(\lambda,\nu)$ to  $((\log|s|)^{-1}\lambda,\log|s|\nu)$ we will
  need to rescale all the parameters as well:
  \begin{gather*}
    \lambda\sim\lambda_0 \log|s|,\ \beta\sim\beta_0 \log|s|,\
    \gamma\sim\gamma_0 \log|s|,\ h\sim h_0\log|s|,\ c\sim c_0
    \log|s|,\\
    \nu\sim\frac{\nu_0}{\log|s|}, \
    \beta^\vee\sim\frac{\beta^\vee_0}{\log|s|}, \
    \Hess\Phi\sim\frac{\Hess\Phi_0}{(\log|s|)^2}.
  \end{gather*}
  Or equivalently, we could apply the log map with the base $|s|$ as
  in \cite{HZh}.

  We saw in the proof of Lemma~\ref{lemma:scalar} that $\Hess
  \Phi^\sm_{\lambda_0+\gamma_0}$ is degenerate along $\F$ outside
  $\Dv_{\lambda_0+\gamma_0+h_0}$. This argument extended to the entire
  toric variety shows that up to terms of order $O(\epsilon)$, the set
  $X_T\setminus\log^{-1} (\Dv_{\lambda_0+\gamma_0+h_0})$ (which
  contains $Z_s$) has distance from $\dDv_{\lambda}$ bounded by the
  diameters of the torus fibers $\T$.
  The size of the tori $\T$ is determined by the norm of
  $\Hess\Phi^\sm_{\lambda+\gamma}$ at the corresponding point, which
  is bounded by $\frac{C_1(\gamma_0)}{(\log|s|)^2}$.

  The rest of the proof consists of careful picks for asymptotics of
  the rescaled parameters
  $\beta_0,\beta^\vee_0,h,\gamma_0,c_0,\epsilon$ to ensure that the
  following expressions
  \begin{gather*}
    \beta_0,\quad \beta^\vee_0,\quad \gamma_0,\quad c_0, \quad h_0,\\
    C_1(\gamma_0)(\log|s|)^{-2},\\
    \log|s|\frac1{\beta^\vee_0} e^{-\beta_0 \log|s|}\cdot C_0
    (\beta_0),\\ 
    \log|s|\frac{ C_0 (\beta_0)C_1(\gamma_0)}{\beta^\vee_0}
    e^{-\beta_0 \log|s|}+ C_2(\gamma_0)+C_3(h_0)\\ 
    \hskip 2in
    +C_4(\gamma_0)+C_5(\beta_0,c_0)+O(\epsilon)(\log|s|)^{2} 
  \end{gather*}
  go to 0 as $\log|s|\to\infty$, where the first two lines take care
  of the Hausdorff convergence, and the last two give matching of the
  complex structure and the metric under the embedding
  $\psi:Z^\sm_s\to W_s$ 
  
  For instance, we can choose
  \begin{equation*}
    \beta^\vee_0\sim \frac1{\log|s|},\quad 
    \gamma_0\sim C_1^{-1}(\log|s|),\quad c_0\sim \frac1{\log|s|},
    \quad h_0\sim \frac1{\log|s|}.
  \end{equation*}
  And $\beta_0(\log|s|)$ has satisfy
  $C_1(\frac1{\log|s|}e^{-\beta_0\log|s|})<\log|s|$
  (i.e. $e^{-\beta}<\gamma$, which is needed for the proof of
  Lemma~\ref{lemma:scalar}) and
  $$(\log|s|)^3C_0(\beta_0) e^{-\beta_0\log|s|}\to 0,$$
  which is possible due to the fast decreasing factor of
  $e^{-\beta_0\log|s|}$ when $\beta_0(\log|s|)$ is changing slowly.

Finally, notice that the bi-PIKAS of type $((\log|s|)^{-1}\cdot\lambda, \log|s|\cdot\nu)=(\lambda_0+(\log|s|)^{-1}\log|a|,\nu_0)$ converges to the $(\lambda_0,\nu_0)$-bi-PIKAS.
\end{proof}

\begin{remark}
  We can rephrase the above theorem in terms of the alternate
  definition of the torus bundles
  $W_{\frac1{\log|s|}}(\lambda_0,\nu_0)$ associated to the given
  K\"ahler affine structure on $\Sigma$. Then the statement of the
  theorem will coincide with the Conjecture~2 of \cite{KS}.
\end{remark}

\bibliographystyle{alpha}
\bibliography{references}

\end{document}